\documentclass[12pt]{article}

\usepackage[top=3.7cm, bottom=3.7cm, left=2.8cm, right=2.4cm]{geometry}

\usepackage{setspace}
\usepackage{footmisc}
\makeatletter
\newcommand\iraggedright{%
\let\\\@centercr\@rightskip\@flushglue \rightskip\@rightskip
\leftskip\z@skip}
\makeatother
\iraggedright

\title{Why Did Weyl Think that Emmy Noether Made Algebra the Eldorado of Axiomatics?}

\author{Iulian D. Toader\thanks{\bigskip Institute Vienna Circle, University of Vienna, iulian.danut.toader@univie.ac.at}}

\date{}

\begin{document}

\maketitle

\begin{quote}
    
\textbf{Abstract}: The paper attempts to clarify Weyl's metaphorical description of Emmy Noether's algebra as the Eldorado of axiomatics. It discusses Weyl's early view on axiomatics, which is part of his criticism of Dedekind and Hilbert, as motivated by Weyl's acquiescence to a phenomenological epistemology of correctness, then it describes Noether's work in algebra, emphasizing in particular its ancestral relation to Dedekind's and Hilbert's works, as well as her mathematical methods, characterized by non-elementary reasoning, i.e., reasoning detached from mathematical objects. The paper turns then to Weyl's remarks on Noether's work, and argues against assimilating her use of the axiomatic method in algebra to his late view on axiomatics, on the ground of the latter's resistance to Noether's principle of detachment.
\end{quote}

\section{Introduction}

Weyl's memorial address for Emmy Noether contains the following description of the axiomatic method in mathematics: ``Hence axiomatics is today by no means merely a method for logical clarification and deepening of the foundations, but it has become a powerful weapon of concrete mathematical research itself. This method was applied by Emmy Noether with masterly skill, it suited her nature, and she made algebra the Eldorado of axiomatics.'' (Weyl 1935, 60) On the face of it, this metaphor is not unequivocal. Did Weyl mean that Noether transformed algebra into a land of riches, a domain where great mathematical wealth can be rapidly and rather effortlessly acquired by means of axiomatics, or did he mean that she turned it into something like a mythical realm that keeps luring mathematicians with the illusory promise of bountiful fortunes? Moreover, what motivated this metaphorical description in the first place? 

In a summer course that took place in October 1931, Weyl presented a paper titled ``Topology and Abstract Algebra as Two Roads of Mathematical Understanding" to high-school teachers in Switzerland, with whom he shared the following thought: ``today the feeling among mathematicians is beginning to spread that the fertility of these abstracting methods [such as the axiomatic method] is approaching exhaustion. ... I foresee that the generation now rising will have a hard time in mathematics.'' (Weyl 1932, 47; quoted from Weyl's own translation in 1935, 60sq.) This prediction surely did not sound too encouraging to at least some of the teachers in that audience. Weyl quite bluntly criticized the axiomatic method for having no independent ability to produce results, for having to rely on the substantial wealth acquired by non-axiomatic methods, as well as for draining this mathematical wealth. The alarming signal to mathematicians was clear enough: avoid axiomatics as far as possible.

A while after Weyl's paper has been published, Emmy Noether expressed her disagreement with Weyl's view on axiomatics, quite plausibly in July 1933 at \textit{Gasthof Vollbrecht} near G\"ottingen, where they met with the occasion of Emil Artin's visit (see Eckes and Schappacher 2016, for the exact dating of a famous photo taken with that occassion). A couple of years later, after they had already emigrated to the United States, Weyl duly recalled their discussion: ``Emmy Noether protested against that [i.e., against his contention that axiomatics is approaching exhaustion]: and indeed she could point to the fact that just during the last years the axiomatic method had disclosed in her hands new, concrete, profound problems ... and had shown the way to their solution.'' (Weyl 1935, 61). In all appearance, Noether's protest made Weyl reconsider his view about the axiomatic method. But why exactly did it make him do so, and to what extent did his view actually change? Since there is no direct textual evidence of what Noether said to Weyl, this apparent reconsideration must be reconstructed, for the most part, from his reflections on axiomatics, both before and after their \textit{Ausseinandersetzung}. 

As is quite well known, in the first part of his career Weyl criticized the views of mathematics he attributed to Dedekind and Hilbert. This criticism, which I discuss in sections 2 and 3, discloses a rather negative disposition towards modern axiomatics. More specifically, Weyl criticized what he considered to be Dedekind's epistemologically perverted view on proof, characterized by a reversal of what Weyl took to be the right epistemology of mathematical belief. He also criticized his own \textit{Doktorvater}, Hilbert, for the limited epistemological traction of his axiomatic method. As we will see, it was Weyl's acquiescence to a phenomenological epistemology of mathematics that played a major role in motivating both these criticisms, and made him express a considered preference for ``elementary'' reasoning in mathematics, i.e., reasoning ``undetached'' from objects. This preference, as we will further see, led him to a strongly revisionist attitude towards algebra.

In section 4, I discuss Emmy Noether's mathematical methods, as especially illustrated by her work in invariant theory and algebra, and I note precisely the connection to Dedekind and Hilbert. Noether, herself, often emphasized her indebtedness to their mathematical work. ``Es steht schon bei Dedekind'', she used to say, apparently (see, e.g., Corry 2017, 134). Philosophical discussions of her algebraic methods are scarce, but those few provided by some of Noether's collaborators and by contemporary commentators converge on the point that what characterizes these methods is precisely the extensive use of what Weyl would call non-elementary reasoning, i.e., reasoning detached from mathematical objects. Indeed, I suggest that, for understanding Noether's discussion with Weyl and the subsequent change in his thinking about axiomatics, this is the most salient epistemological feature of her methods that needs to be appreciated. I further present Weyl's characterizations of Noether's axiomatics, in his funeral speech, on April 17, 1935, and in his memorial address delivered a week later. On both occasions, he emphasized the change her methods effected in algebra and the lasting and growing influence he thought they would continue to have, while at the same time pointing out some of what he considered to be their shortcomings. But it is hard to see how Noether's methods could have led Weyl to reconsider his take on axiomatics, given the differences between their relations to the views of Dedekind and Hilbert.

In section 5, I turn to Weyl's later work, which shows that his early negative disposition towards axiomatics changed indeed, but it changed in the following way: while continuing to maintain a critical attitude with respect to Hilbert's axiomatic method, Weyl came to distinguish it from the type of axiomatics that he saw at work in algebra. I argue that, despite his own pronouncement, the latter type of axiomatics cannot be identified with Noether's axiomatics. For the most characteristic feature of Noether's axiomatics -- detachment of reasoning from objects and states of affairs -- was never acceptable to Weyl. He thought this feature must be eliminated if genuine mathematical knowledge is to be obtained. This, I contend, throws light on his metaphorical description of abstract algebra as the Eldorado of axiomatics: for Weyl, this turns out to be a domain where great mathematical wealth could truly be earned, provided that the non-elementary character of reasoning can be eliminated. But if this requirement is unsatisfied, the Eldorado remains a mere myth that continues to seduce and delude mathematicians into believing that great wealth is literally within their reach.

\section{Weyl's early view on axiomatics (I)}

Let us start with an important passage from a footnote in \textit{Das Kontinuum}: ``In the Preface to the first edition of Dedekind's famous \textit{Was sind und was sollen die Zahlen?}, we read that `In science, what is provable ought not to be believed without proof.' This remark is certainly characteristic of the way most mathematicians think. Nevertheless, it is a perverse principle. As if such a indirect concatenation of grounds as what we call a `proof,' can awaken any `belief' without our assuring ourselves, through immediate insight, of the correctness of each individual step! This (and not the proof) remains throughout the ultimate source of knowledge; it is the experience of truth.'' (Weyl 1918, 119; translation amended) What Weyl seems to mean here is that what one might call a Dedekindian proof, i.e., a logical derivation of a mathematical theorem from the basic laws of thought, cannot make us properly believe a theorem. For the merely logical derivation from such grounds would not be enough to assure us of its epistemological correctness. The proper belief of a mathematical theorem, Weyl emphasized, requires an epistemologically correct proof. In contrast, according to what he took to be Dedekind's norm of mathematical belief, proper belief requires Dedekindian proof.\footnote{For a detailed analysis of Weyl's footnote, see Toader 2016. Translating Weyl's ``ein verkehrtes Prinzip'' by ``a perverse principle'' captures better his intended meaning than the 1987 English translation --- ``a preposterous principle''. I think Weyl did not mean that Dedekind's principle was simply absurd, that it did not make sense. Rather, he meant that it was contrary to sense, or backwards, that it expressed the reverse of (what he thought was) the right view on mathematical belief. An alternative translation, suggested by one referee, could be ``an inverted principle''.}

How does this phenomenological requirement of epistemological correctness reflect on the nature of mathematical reasoning? In the main text of \textit{Das Kontinuum}, Weyl presented his normative view on proof as a series of correct judgments, in the following terms: ``The presentation of the fact that a judgment U is a consequence of the axioms can and must be done ... through a generally ramified organism of `elementary' inferences, which in order to be communicated must be artificially transformed in a chain of interlocking links. It is in this way that mathematical \textit{proof} comes about; in it, all insight that is to be satisfied is concentrated on the logical inferences and is no longer directed at the objects and states of affairs judged upon.'' (\textit{ibid.}, 17; translation amended) Weyl deplored the fact that proof is often taken to be the result of a process of detachment from mathematical objects and states of affairs. Such detachment renders the self-evidence of a mathematical judgment inaccessible to intuition or immediate insight, and therefore prevents one's having an experience of its truth. On Weyl's view, genuine proof requires elementary reasoning, i.e., reasoning that is undetached from mathematical objects and states of affairs, for only such reasoning can constitute an epistemologically correct proof. More exactly, the requirement is that reasoning ought to be conducted within a fully interpreted language, in which all mathematical symbols have meaning, thereby allowing for the possibility of expressing self-evident judgments.\footnote{One might immediately ask whether predicative proofs meet the standards imposed by the phenomenological epistemology advocated by Weyl. While this issue can only be briefly addressed  below, it seems to me that this epistemology would invalidate much more of analysis than predicativists, including Weyl himself, would be willing to give up.}

Where is the perversity, one might ask. Weyl considered  epistemologically perverse or backwards the requirement that, if proper belief of theorems is to be achieved, Dedekindian proofs should replace (what Weyl took to be) epistemologically correct proofs. But he thought that this replacement stemmed from a mistaken epistemology of mathematical belief. What led him to think in this way was, in all appearance, his acquiescence to a phenomenological epistemology of mathematics. The most important element of this epistemology, which originates with Brentano, and the only one that we need to emphasize here is correctness (\textit{Richtigkeit}), which is primarily taken to be a property of judgments. A judgment is correct if and only if it possesses self-evidence (\textit{Evidenz}), which is epistemically accessible in intuition or immediate insight, and it is the core of what Weyl calls an experience of truth, using another phenomenological term appropriated from Husserl.\footnote{For a more detailed analysis of Weyl's involvement with Husserl's epistemology, see Toader 2011. For Brentano's notion of \textit{Richtigkeit}, see, e.g., Textor 2019.} According to this phenomenological epistemology, a proof should provide an experience of truth, and this can occur if and only if the proof is correct, i.e., a series of (all and only) correct judgments. Dedekindian proof, Weyl argued, not only falls short with respect to this phenomenological norm of belief, but even more than that, Dedekind's principle advocates what Weyl considered to be a reversal of that norm. Hence, his calling it epistemologically perverse, i.e., not merely senseless, but contrary to sense.

Weyl's project in \textit{Das Kontinuum}, as is well known, was the reconstruction of analysis in the face of impredicativity, one in which he attempted to ``replace this shifting foundation [of analysis] with pillars of enduring strength'' (Weyl 1918, 1). As he emphasized in another fitting metaphor, this was in fact an amputation meant to excise ``every cell (so to speak) of this mighty organism [of analysis] permeated by the poison of contradiction'' (\textit{ibid.}, 32). Less metaphorically, the project aimed at providing a justification to the claim that Dedekind's definition of real numbers as sets of rational numbers is possible only if certain restrictions are imposed on what is to count as a definition. Analysis is to start from a basic category of mathematical objects -- the totality of natural numbers -- and to extend this category via explicit definitions. A fundamental feature of this extension is that it excludes impredicative definitions, i.e., definitions that make reference to sets of objects in which the defined object is an element. This restriction entails that some theorems of classical analysis, like the least upper bound principle, cannot be proved for arbitrary bounded sets of real numbers. However, many theorems of classical analysis can be proved for predicatively definable bounded sets of real numbers, that is for countable sequences of reals: for example, the least upper bound principle, but also the entire classical theory of continuous functions, and based on this, the fundamental theorem of algebra. 

Logicians have long focused on the predicative restrictions on classical analysis and, thus, on the extent of Weyl's predicative analysis.\footnote{See, e.g., Grzegorczyk 1955 and Feferman 1964.} In a letter to G\"odel, from 1932/1933, prompted by the publication of the second edition of Weyl's book, Karl Menger expressed a characteristic complaint: ``If only you could help me to formulate the axioms so that the set of all endpoints of a curve is admissible!! Dear Mr. G\"odel, do ponder these (for me) terribly tormenting matters! And don't laugh at me! No one besides me is more knowledgable about it. You saw, I'm sure, that Fraenkel didn't know that one can extract the entire continuum without the axiom of reducibility. So even if it looks like a triviality to you, do take the trouble to do this. In the end, the matter must be clarified for mathematicians who are not specialists in logic. Just think of Weyl's disgusting prattle  (\textit{widerliches Geschw\"atz})!'' (G\"odel 2014, 101) Disgusting or not, Weyl's predicative system allows no type-theoretical reconstruction of analysis and includes no Russellian axiom of reducibility. It allows no unrestricted second-order quantifiers, and so it is typically thought that its underlying logic could be either predicative second-order logic or classical first-order logic. Most commentators tend to believe that Weyl adopted the latter.\footnote{See, e.g., Feferman 1988.} 

To be sure, Weyl's system is second-order, one in which natural numbers are first-order objects. Second-order quantifiers, to the extent that they are allowed, must be restricted to the predicatively definable subsets of natural numbers. But since the totality of natural numbers is assumed to be given, the system can be taken to validate the $\omega$-rule for first-order quantifiers, the same way in which it validates the arithmetical axioms, including complete induction. For any formula $A(x)$ in the language of arithmetic, this rule licenses the inference from  $A(0), A(1), A(2), ...,$ to $\forall x A(x)$. The logic of the system is, thus, what might be called $\omega$-logic.\footnote{See Moore 1988, 135. For the history of the $\omega$-rule, see e.g. Buldt (2004).} But if that is the case, then arguably no predicative proof in the system would be able to provide an experience of truth, since the $\omega$-rule is infinitary (i.e., any inference based on it has an infinite number of premises). This assumes, of course, that the application of an infinitary rule entails detachment from mathematical objects and states of affairs, which as already noted makes self-evidence inaccessible. In Weyl's own terms, it assumes that reasoning in accordance with the $\omega$-rule is not elementary. If the assumption is true, as seems to be the case, then predicative proofs support theorems that exceed the boundary conditions on proper mathematical belief imposed by Weyl's own phenomenological epistemology. As a consequence, one would have to give up these theorems, including the fundamental theorem of algebra. But that would seem to be a rather excessively revisionist attitude, in the sense that very few mathematicians would want to adopt it. Nevertheless, as we will see presently, Weyl would as a matter of fact come to argue, in the early 1930s, that elementariness requires that the fundamental theorem of algebra should indeed be given up.

\section{Weyl's early view on axiomatics (II)}

The phenomenological epistemology of mathematics was the basis not only for Weyl's criticism of Dedekind's conception of proof, but also for his later emphasis on the limited epistemological traction of Hilbert's axiomatic method. Weyl's often quoted complaint about this method is as follows: ``Hilbert's mathematics may be a pretty game with formulas, more amusing even than chess; but what does it have to do with knowledge, since its formulas should admittedly have no contentual significance by virtue of which they would express intuitive truths?'' (Weyl 1927, 61) The objection here is that the kind of reasoning that characterizes part of Hilbert's mathematics, i.e., non-elementary reasoning within a partially uninterpreted language, in which not all symbols have meaning, lacks the ability to express self-evident judgments. As a consequence, Weyl argued, Hilbert's mathematics cannot provide epistemologically correct proofs and, thus, cannot be a source of mathematical knowledge.

Weyl continued to be critical of Hilbert's views in the philosophy and foundations of mathematics. The full passage from the 1932 paper quoted above, in his own 1935 translation, included the following criticism: 

\begin{quote}
    
I should not pass over in silence the fact that today the feeling among mathematicians is beginning to spread that the fertility of these abstracting methods is approaching exhaustion. The case is this; that all these nice general notions do not fall into our laps by themselves. But definite concrete problems were first conquered in their undivided complexity, single-handed by brute force, so to speak. Only afterwards the axiomaticians came along and stated: Instead of breaking in the door with all your might and bruising your hands, you should have constructed such and such a key of skill, and by it you would have been able to open the door quite smoothly. But they can construct the key only because they are able, after the breaking in was successful, to study the lock from within and without. Before you can generalize, formalize and axiomatize, there must be a mathematical substance. I think that the mathematical substance in the formalizing of which we have trained ourselves during the last decades, becomes gradually exhausted. And so I foresee that the generation now rising will have a hard time in mathematics. (Weyl 1932, 47; 1935, 60sq.)

\end{quote}

Weyl considered that Hilbert's axiomatic method is dependent on the conceptual wealth acquired by non-axiomaticians and has no independent ability to produce mathematical results. He also gave examples to illustrate such criticism. As late as in his speech as the Fields Committee President, for instance, Weyl would note: ``Only if someone has the courage of attacking the primary concrete problems in all their complexity, will the general concepts gradually emerge which resolve the difficulties and ease the further progress.'' (Weyl 1954, 619) The concept that Weyl had in mind as an example is that of a Kodaira dimension, a numerical invariant on the basis of which Kodaira classified algebraic varieties, such that each family of algebraic varieties has a given Kodaira dimension. The emergence of this concept allowed, as Weyl remarked, a ``profound generalization of the well-known fact that every compact Riemann surface belongs to an algebraic function field.'' (\textit{loc. cit.}) The notion of a Kodaira dimension is supposed to illustrate how mathematical concepts are obtained by fertile, non-axiomatic methods, rather than falling into the axiomatician's lap by themselves. 

Weyl's grim prediction from that October 1931 was based on what he saw as a misguided view of concept-formation in mathematics, but the core of his complaint against the axiomatic method seems to have been that this fails to adhere to the epistemological principles of correctness and elementariness, that it allows non-elementary reasoning, i.e., reasoning within partially uninterpreted languages, and consequently that it makes mathematical proof unable to provide an experience of truth. Weyl's commitment to using only elementary reasoning in proofs led him to advocate a revisionist attitude that, as I noted above, very few mathematicians would find tolerable. For this commitment made him think that one would do well to avoid using the fundamental theorem of algebra in proofs: 

\begin{quote}
    
If one operates in an arbitrary abstract number field rather than in the continuum
of complex numbers, then the fundamental theorem of algebra, which asserts
that every complex polynomial in one variable can be [uniquely] decomposed into
linear factors, need not hold. Hence the general prescriptíon in algebraic work:
See if a proof makes use of the fundamental theorem or not. In every algebraic
theory, there is a more elementary part that is independent of the fundamental
theorem, and therefore valid in every field, and a more advanced part for whích
the fundamental theorem is indispensable. The latter part calls for the algebraic
closure of the field. In most cases, the fundamental theorem marks a crucial
split; its use should be avoided as long as possible. (Weyl 1932, 46) 
\end{quote}

Thus, on Weyl's view, elementary reasoning -- in this case, reasoning that is undetached from the elements of a number field -- is always to be preferred to higher or non-elementary reasoning -- in this case, reasoning over the algebraic closure of a field. The assumption here seems to be that idealization procedures like algebraic closure entail detachment from mathematical objects and states of affairs. This assumption was of course justified by Steinitz's 1910 theorem, which states that the algebraic closure of a field is unique up to isomorphism, and which suggests that mathematical reasoning concerns the structure of mathematical relations over an algebraic closure, rather than their relata. It seems plausible that Weyl's insistence on elementary reasoning, undetached from mathematical objects, which led him to adopt what looks like an extremely revisionist attitude in algebra, was here again motivated by his acquiescence to a phenomenological epistemology of mathematics, which had earlier led him to assert the epistemological correctness of proofs as a necessary condition for genuine mathematical knowledge.

Furthermore, it could further be argued that this attitude was, once again, directed against Hilbert. Indeed, Weyl's belief that elementariness overrides any epistemic value the fundamental theorem of algebra might have goes straight against Hilbert's doctrine of ideal elements. According to the latter, the postulation of transfinite axioms is meant to increase the epistemic value of finitary mathematics, e.g., its simplicity and fruitfulness, just like the introduction of complex numbers via the fundamental theorem of algebra was meant to increase the epistemic value of real analysis. As Hilbert famously put it: ``In my proof theory, the transfinite axioms and formulae are adjoined to the finite axioms, just as in the theory of complex variables the imaginary elements are adjoined to the real, and just as in geometry the ideal constructions are adjoined to the actual'' (Hilbert 1923, 1144) The epistemic value of the fundamental theorem of algebra is, of course, typically associated with the simplicity induced by the linear factorization property of polynomials defined over complex fields.\footnote{Cf. Detlefsen 2005. Of course, the linear factorization of a polynomial also obtains on minimal extensions of a number field, such as the splitting field of the polynomial. For a contemporary defense of this Kroneckerian view, according to which the fundamental theorem of algebra (in its standard formulation) is actually false, see Edwards 2005.} But its value can also be associated with the fact established by Steinitz's theorem mentioned above, and its structuralist reading. However, this fact appears to have been regarded by Weyl as an injunction against the elementary character of mathematical reasoning, and motivated his revisionist attitude towards the fundamental theorem of algebra, while strengthening his case against Hilbert's axiomatics and in particular his view on idealization in mathematics. The background to all this, as suggested already, was Weyl's acquiescence to the phenomenological epistemology of mathematics described above in section 2. After 1932, however, Weyl's view on axiomatics changed. In order to see how it changed, and the extent to which it did, we must consider the presumed efficient cause of this change: Emmy Noether's methods in algebra.

\section{Noether's principle of detachment}

In the history of mathematics, Emmy Noether is often referred to as ``the mother of modern algebra'', and thought to have ``changed the face of algebra by her work'', as Weyl himself put it in 1935. The foundations she built upon, however, as she never failed to mention, are two: Dedekind and Hilbert. This is also often emphasized in contemporary reflections on her work: ``One theme associated with the name of Emmy Noether is that of the `abstract axiomatic approach' to algebra; a second is that of attention to \textit{entire algebraic structures and their mappings}, rather than to, say, just numerical attributes of those structures (as a notable example, in topology, attention to homology groups and their induced mappings, not just to Betti numbers and torsion coefficients). One might consider that, of these themes, the second links Noether to Dedekind, the first rather to Hilbert.'' (Stein 1988, 245) There were, of course, other mathematicians who greatly influenced Noether's work, such as Ernst Steinitz, whose paper on the algebraic theory of fields (which includes his categoricity theorem mentioned above) was cited in every mathematical paper that Noether wrote in the decade following its publication (Roquette 2010). But it is fair to say that the work of Dedekind and Hilbert was essential to her own. 

It is well known that Noether's study of invariant theory, for example, had developed from a purely algorithmic approach, such as in her 1907 doctoral dissertation in Erlangen, ``On Complete Systems of Invariants for Ternary Biquadratic Forms'', defended under Paul Albert Gordan, to the conceptual approach illustrated by  Hilbert's first proof of the so-called Finite Basis Theorem, a corollary of which was that every form of any degree and in any number of variables has a finite complete system of invariants. Her work in algebra, as just noted in the quote from Stein, links her to Hilbert's own axiomatic approach. Also well known is that much of Noether's work on commutative algebra, and especially her paper ``Ideal Theory in Rings'' (1921), where she developed her axiomatic theory of ideals, was based on the so-called ascending chain condition, a notion that had been articulated first by Dedekind (in 1894) for rings of algebraic integers. She generalized this condition to arbitrary commutative rings, now called Noetherian rings, and gave it expression in the so-called \textit{Teilerkettensatz}, an axiom that states that any chain of ideals necessarily comes to an end after a finite number of steps if each of them comprises the preceding one as a proper part, i.e., if the chain is ascending. Other notions introduced by Dedekind, which were further developed by Noether, include the concept of a module, which I'll come back to below. As noted by Stein, what links her to Dedekind perhaps most clearly is Noether's emphasis on thinking in terms of mappings between algebraic structures, especially mappings that preserve structural relations. Moreover, her influence on contemporary mathematics should be emphasized as well: consider, for example, the notion of Noetherian topological spaces, i.e., spaces that have the descending chain condition on closed subsets, and its role in the development of the mathematics of Zariski geometries (Zilber 2010).

Now, this is nothing new.\footnote{For a more detailed analysis of Noether's work, see McLarty 2006.} But it is often said that Noether introduced new methods in algebra,  which are characterized in various terms as abstract, axiomatic, structural, and conceptual, and that these methods made algebra so general as to apply to all mathematics. What are these methods? When prompted by her collaborators, she would say things like this: ``Meine Methoden sind Arbeits- und Auffassungsmethoden, und daher anonym \"uberall eingedrungen.''’ (letter to Hasse from 12 Nov. 1931, quoted in Koreuber 2015, 71) Surely, characterizing her methods as \textit{Arbeits- und Auffassungsmethoden} does not fully clarify their nature, but rather indicates why they came to be applied so widely in modern mathematics. More helpful in this respect, van der Waerden formulated in his \textit{Nachruf} for Noether a maxim that, according to him, always guided her work as a central principle for her mathematical thinking: ``All relations between numbers, functions and operations become perspicuous, capable of generalization, and truly fruitful, when they are \textit{detached from their particular objects} and reduced to general conceptual connections.''\footnote{The German original is: ``Alle Beziehungen zwischen Zahlen, Funktionen und
Operationen werden erst dann durchsichtig, verallgemeinerungsf\"ahig und wirklich
fruchtbar, wenn sie \textit{von ihren besonderen Objekten losgel\"ost} und auf allgemeine
begriffliche Zusammenh\"ange zur\"uckgef\"uhrt sind.'' (van der Waerden 1935, 469; my emphasis)} Similar claims have been articulated in more recent characterizations of Noether's work, such as: ``to get abstract algebra away from thinking about operations on elements, such as addition or multiplication of elements in groups or rings. Her algebra would describe structures in terms of selected subsets (such as normal subgroups of groups) and homomorphisms.'' (McLarty 2006, 188) All this arguably points in the direction of a certain structuralist epistemology of mathematics, one that considers relations, rather than their relata, as the main focus of mathematical investigation and the main source of mathematical knowledge. However, what seems to me to be the real significance of these descriptions of Noether's methods is their emphasis on the epistemic value of detachment from mathematical objects, from the particular elements of number fields, groups, rings, etc., and from the operations on them. Such detachment is considered as the main source for the fruitfulness, perspicuity, generality, and simplicity of Noether's mathematical results in algebra. It is this -- what I want to call the principle of detachment -- that I take to characterize Noether's methods as applied in her work in algebra, and that I think needs to be appreciated in the context of her \textit{Ausseinandersetzung} with Weyl. If this is true, it's then easy to imagine the proportions of Noether's dissatisfaction with a view like the one propagated by Weyl, which she had seen most recently advocated in his 1932 paper -- the view that detachment from mathematical objects is, on the contrary, to be avoided precisely because of its epistemic shortcomings. 

In his funeral speech and in the subsequent memorial address, Weyl recorded for the first time Noether's opposition to his view. He recalled his previous emphasis on the epistemological correctness of proofs and the elementary character of reasoning, as well as his criticism of the axiomatic method, stressing again the fruitlessness and dependence of this method on non-axiomatic methods. Then he noted the following: ``Emmy Noether protested against that: and indeed she could point to the fact that just during the last years the axiomatic method had disclosed in her hands new, concrete, profound problems by the application of non-commutative algebra upon commutative fields and their number theory, and had shown the way to their solution.'' (Weyl 1935, 61) One can imagine that Noether had pointed out that her focus on abstract conceptual structures, motivated by the principle of detachment, not only led to a better understanding of old mathematical problems and their solutions, but also to uncovering new ones. Thus, it was not exhaustion, but rather fruition, that her axiomatics was coming to.

Noether's axiomatic method was illustrated by Weyl with reference to her 1920 paper, ``On Modules in Non-Commutative Fields, Especially from Differential and Difference Expressions'' (co-authored with Werner Schmeidler): 

\begin{quote}
    
It is here for the first time that the Emmy Noether appears whom we all know, and who changed the face of algebra by her work. Above all, her conceptual axiomatic way of thinking in algebra becomes first noticeable in this paper dealing with differential operators as they are quite common nowadays in quantum mechanics. In performing them, one after the other, their composition, which may be interpreted as a kind of multiplication, is not commutative. But instead of operating with the formal expressions, the simple properties of the operations of addition and multiplication to which they lend themselves are formulated as axioms at the beginning of the investigation, and these axioms then form the basis of all further reasoning. A similar procedure has remained typical for Emmy Noether from then on. (Weyl 1935, 55sq.) 

\end{quote}

Noether's paper starts from the abstract definition of an arbitrary field, and then on the basis of the non-commutative character of multiplication, defines various types of modules and proves that two modules are of the same type if and only if their sets of residue classes are isomorphic, i.e., if the sum of any two residue classes from the first set maps to the sum of the corresponding residue classes from the second set, and the product of any residue class from the first set with a polynomial maps to the product of the corresponding residue class from the second set with the same polynomial. Various other results are then proved in the general theory of modules. 

But why would such results, and others that were to follow the same method, convince Weyl to reconsider his view of axiomatics? The answer lies, I believe, in his emphasis on what he took to characterize Noether's axiomatic method: its focus is not on algebraic operations on field elements, but rather on the algebraic relations between these elements and on the classes of abstract structures (and their isomorphism types) that can be defined thereupon. Weyl was able to see here Noether's principle of detachment at work. He further noted in his funeral speech: ``No-one, perhaps, contributed as much as [Emmy Noether] towards remoulding the axiomatic approach into a powerful research instrument, instead of a mere aid in the logical elucidation of the foundations of mathematics, as it had previously been.''\footnote{See Roquette 2008. This includes Weyl's entire funeral speech, which was preserved in Grete Hermann's archive, in the ``Archiv der sozialen Demokratie'' in Bonn.} Thus, Weyl seems to have been able to recognize the strength of detachment and the epistemic value of non-elementary reasoning in mathematics, suggesting here for the first time what he would later articulate as an important distinction between Hilbert's and Noether's axiomatic methods. A week later, in the memorial address, Weyl took time to emphasize again this distinction, and then described Noether's algebra as the Eldorado of axiomatics. As mentioned at the outset, this metaphorical description is not unequivocal: does it mean that, thanks to axiomatics, a great wealth of mathematical results can now be rapidly and easily obtained in algebra, or does it mean that, despite axiomatics, all that wealth is bound to remain a tortuous delusion, a mere myth? 

To be sure, the same metaphor had been included in Weyl's 1932 paper, without explicitly mentioning Emmy Noether as the mathematician most responsible for this transformation of algebra. After comparing the visual simplicity and intelligibility of topological methods with the algorithmic methods employed in algebra, he noted: ``Now,
step by step, abstract algebra tidies up the clumsy computational apparatus. The generality of the assumptions and axiomatization force one to abandon the path of blind computation and to break the complex state of affairs into simple parts that can be handled by means of simple reasoning. Thus algebra turns out to be the El Dorado of axiomatics.'' (Weyl 1932, 40) But, again, it's not clear what Weyl took this to mean. For in the same paper, as we have seen, he criticized axiomatics for its lack of any independent ability to produce results, and for its inevitable dependence on concepts introduced non-axiomatically. He took fruitfulness as ``perhaps the only criterion'' that could justify generalization and formalization, but warned that the fertility of these procedures that characterize axiomatics lead to exhaustion, rather than a wealth of mathematical results. Furthermore, while his topological example of fruitful theories was Riemann's theory of algebraic functions and their integrals, Weyl's algebraic example was Dedekind's theory of ideals, introduced in the latter's theory of abstract number fields and subsequently developed by Emmy Noether. A close reading, however, indicates that Weyl took the latter development to run into what he then considered to be the problem of detachment, whose solution he thought was the general prescription that we discussed above: avoid using the fundamental theorem of algebra in proofs as long as possible, because reasoning that involves algebraic closure is non-elementary reasoning, i.e., the kind of reasoning that Weyl thought ought to be resisted in mathematics. 

Are Weyl's later pronouncements on axiomatics significantly different? Did he really come to fully appreciate and embrace Noether's principle of detachment? I think the answer is no, as I will try to argue in the rest of the paper. Even in the 1935 memorial address, Weyl promptly noted that Noether's ``strong drive toward axiomatic purity'' was not without shortcomings, for the application of her principle of detachment and the abstractness of her concepts ``had the disadvantage ... that she was sometimes but incompletely cognizant of the specific details of the more interesting applications of her general theories.'' (Weyl 1935, 63) Then he went again metaphorical: ``Her methods need not, however, be considered the only means of salvation.'' Whenever available, he suggested, constructive methods are preferable to axiomatic ones due to the elementary character of constructive reasoning, i.e., to its being undetached from mathematical objects and states of affairs. Thus, it seems really hard to see how Weyl, who had so forcefully criticized both Dedekind and Hilbert, would come to fully appreciate Noether's results and methods, based as they were on results and methods initiated and developed precisely by Dedekind and Hilbert. The fact that he seems to have done so cries out for an explanation.

\section{Weyl's later view on axiomatics}

In a later paper, Weyl came back to his distinction between the different roles of the axiomatic method in mathematics: ``One very conspicuous aspect of twentieth century mathematics is the enormously increased role which the axiomatic approach plays. Whereas the axiomatic method was formerly used merely for the purpose of elucidating the foundations on which we build, it has now become a tool for concrete mathematical research. It is perhaps in algebra that it has scored its greatest successes.'' (Weyl 1951, 159) Such passages suggest that Weyl turned from complaining about exhaustion via the axiomatic method to valuing the fruition to which this method led to in algebra. This would be, without doubt, a significant turn, especially if it entailed that what he took to be the greatest successes of axiomatics in Noether's algebra were mathematical results that he, himself, would have also believed as epistemologically correct. But to think so would fly in the face of Weyl's earlier criticism of Dedekind and Hilbert, would ignore his views on proof, on correctness and elementariness, and forget about his acquiescence to a phenomenological epistemology. Of course, it may be conjectured that as his philosophical thought on modern science and mathematics matured, Weyl came to see the shortcomings of that epistemology, to realize that modern mathematics is characterized by non-elementary reasoning, concerned with abstract structures and relations, rather than by elementary reasoning, undetached from their relata. It may also be thought that Weyl came to see Dedekind's principle as expressing a new norm in mathematics, and to realize that his general prescription about the fundamental theorem of algebra was excessively revisionist. But these conjectures are unsupported by available evidence.\footnote{For a discussion of what Weyl took to be the cognitive significance of what he called non-elementary mathematics, when this is taken as a part of a theoretical whole that also includes modern physics, see Toader 2013 and 2014.} In order to see what Weyl's later view on axiomatics actually was, I want to look more closely at a paper published posthumously, but written in the early 1950s: ``Axiomatics Versus Constructive Procedures in Mathematics''. 

In this paper, Weyl introduced a distinction between two uses of the  axiomatic method in mathematics, or two types of axiomatics defined by these uses, which he suggestively called ``transcendental' and ``immanent'' axiomatics. He clarified this distinction by contrasting the epistemological and metatheoretical  properties of two types of axioms. Most importantly, with respect to the method he associated with Hilbert, he noted: ``Transcendental axiomatics has something paradoxical and shocking, because one must try to learn to abstract radically from the familiar intuitive significance of the terms occurring in the axioms as undefined concepts.'' (Weyl 1953, 195) What appeared to him paradoxical and shocking was the fact that, in a Hilbertian system, the axioms are considered as partially uninterpreted sentences, and the proofs that assume them as hypotheses are essentially characterized by detachment from mathematical objects and states of affairs. In Weyl's own terms, then, what was paradoxical and shocking about transcendental axiomatics was the extensive use of non-elementary reasoning, constrained merely by metatheoretical requirements such as consistency, independence, and completeness. In contrast to this, ``immanent axiomatics as applied in algebra or topology or some other concrete branch of mathematics is neither based on external evidence nor on hypotheses, but the axioms are proved to hold for the mathematical objects in the individual situations to which the axioms are applied.'' (\textit{loc cit.}) Unlike transcendental axioms, the immanent ones are not metatheoretically constrained stipulations or hypotheses, and they are not empirical statements either, but they admit of proof nonetheless. The proof of an immanent axiom, I take Weyl to suggest here, is meant to verify its validity in a certain domain of objects, or what he earlier called its epistemological correctness. An immanent axiom is thus taken to be a fully interpreted statement, undetached from mathematical objects and the operations on them. In Weyl's earlier phenomenological terminology, the proof of an immanent axiom would be said to constitute an experience of its truth. Hence, he added: ``The purpose of this sort of axiomatics is to \textit{understand}.'' (\textit{loc cit.}) Immanent axiomatics is meant to provide understanding, one could say, for it is designed to fully satisfy the requirement of elementariness of reasoning, without which no genuine mathematical knowledge or understanding would be possible.

The two types of axiomatics are also distinguished from a metatheoretical point of view: ``While in the [transcendental] axiomatics one was mostly concerned with axioms which determine the structure of the system completely as, e.g., the axioms of Euclidean geometry do for Euclidean space, we have here, in algebra, to do with [immanent] axioms satisfied by many different individual number fields that are not mutually isomorphic.'' (\textit{ibid.}, 194) While completeness, here conceived of as categoricity, is considered an important property of a Hilbertian system, just like consistency and independence, whose importance justifies the metatheoretical anxiety of the transcendental axiomatizer, Weyl underlined the metatheoretical nonchallance of the immanent axiomatizer: ``The questions of consistency, independence and completeness are here [in immanent axiomatics] but of minor importance.'' (\textit{loc. cit.}) Since immanent axiomatics aims at understanding, to say that consistency, independence, and categoricity are of minor importance is to say, at least in part, that they are of minor importance with respect to understanding. This might be taken to mean that consistency, independence, and categoricity are unnecessary for understanding, which entails that understanding can be obtained despite the inconsistency, the lack of independence, or the non-categoricity of a system of immanent axioms. It is arguably hard to see how any real understanding can be obtained on the basis of an inconsistent system of axioms. But Weyl's claim that metatheoretical properties are of minor importance might be also taken to mean that it is proofs of consistency, independence, and categoricity that are unnecessary for understanding. It's not obvious that mathematical understanding should be blocked by the absence of a consistency proof. However, a categoricity proof would apparently constitute an obstacle to understanding, for just like the uniqueness (up to isomorphism) of the algebraic closure of a field, the categoricity of a system of axioms entails detachment from objects, which renders mathematical reasoning non-elementary. Once a categoricity proof is given, attention is arguably redirected away from the objects to the structure of mathematical relations between them. This is why Weyl insisted that immanent axiomatics in algebra is concerned with the individual fields of rational numbers, the real numbers, etc., which are not isomorphic to one another, rather than with isomorphic algebraically closed fields. So it appears that Weyl's revisionist attitude towards algebra, which we recall from his 1932 paper, is preserved in his late view on immanent axiomatics. The immanent axiomatizer follows Weyl's early prescription to avoid the fundamental theorem of algebra, and more generally, to avoid detachment of reasoning from objects. For otherwise, immanent proofs would fall short of epistemological correctness, which would ultimately make immanent axiomatics fail in attaining its very purpose, that of mathematical understanding.

It would be tempting to follow Weyl in thinking that immanent axiomatics, as just described, is nothing other than Noether's axiomatics: ``the immanent axiomatic approach to algebra was developed ... by Emmy Noether and her school.'' (\textit{ibid.}, 198) But doing so would be quite puzzling. As we have seen in section 4 above, Noether's own students and collaborators (as well as contemporary commentators) converge on the view that her axiomatics is essentially characterized by the principle of detachment. This is what turned it, according to Weyl's own assessment, into a powerful research instrument in mathematical practice. In the paper under discussion, however, he associated no lesser epistemic value than mathematical understanding with immanent axiomatics because of its rejection of the principle of detachment, indicated by his characterization of the epistemological and metatheoretical properties of immanent axioms. If this is right, then Weyl's immanent axiomatics should not to be identified with Noether's axiomatics. 

In order to understand the relation between Weyl's immanent axiomatics and Noether's axiomatic approach in abstract algebra, I think one must distinguish between different attitudes towards the principle of detachment. Whereas it is fairly clear that Noether's approach is characterized by detachment as a principle of theoretical rationality, i.e., one without which abstract algebra could not properly develop as a genuine science in its own right, this conception of detachment was never fully acceptable to Weyl. He conceded detachment in the following words: ``We now come to the decisive step of mathematical abstraction: we forget about what the symbols stand for. The mathematician is concerned with the catalogue alone; he is like the man in the catalogue room who does not care what books or pieces of an intuitively given manifold the symbols of his catalogue denote. He need not be idle; there are many operations which he may carry out with these symbols, without ever having to look at the things they stand for.'' (Weyl 1940, 76sq.) Weyl accepted that mathematicians often use partially uninterpreted languages, and realized that sometimes the mathematician does walk ``the path of blind computation''. This clearly indicates that he did come to acknowledge that mathematical reasoning sometimes has a non-elementary character, that it sometimes proceeds by detachment from objects. Nevertheless, he never accepted the view that this is how mathematical reasoning \textit{should} proceed, if its goal is mathematical understanding. For as we have seen, no mathematical understanding could be obtained except by means of elementary reasoning and epistemologically correct proofs. Thus, I think that Weyl, unlike Noether, ultimately regarded detachment as an obstacle that must be eliminated on account of its epistemological shortcomings. His final account was this: ``large parts of modern mathematical research
are based on a dexterous blending of constructive and axiomatic procedures. One
should be content to note their mutual interlocking. But temptation is great,
and not all authors have resisted it, to adopt one of these two views as the genuine
primordial way of mathematical thinking to which the other merely plays a subservient role. ... my own heart draws me to the
side of constructivism.'' (Weyl 1953, 202). In these terms, the difference I see between Noether's axiomatics and immanent axiomatics is the following: the latter is the result of Weyl's having the former play a subservient role to constructivism. This allows Noether's axiomatics as a powerful research instrument in abstract algebra, but denies it the ability to comply with the epistemological strictures underscored by Weyl's unabated constructivism. When these strictures are enforced, Noether's axiomatics becomes immanent axiomatics.

Characterized in this way, immanent axiomatics turns out to be quite similar to the kind of axiomatics that Weyl had described already in \textit{Das Kontinuum}: its axioms are epistemologically correct judgments, and the proofs based on them are what he referred to as ``organisms'', i.e., inferences characterized by elementariness. Such proofs, and only such proofs, are capable of providing an experience of truth: they are the carriers of mathematical understanding. This is how Weyl's acquiescence to phenomenological epistemology made him conceive of axiomatics in 1918. But it seems to me that the same conception of axiomatics, or one quite similar to it, is advocated for in his last paper on this topic. Weyl's immanent axiomatics is, therefore, neither the modern axiomatics of Hilbert, nor that of Noether. Furthermore, I think that clarifying its epistemological character helps disambiguate the description of algebra as the Eldorado of axiomatics, in the following sense: if Weyl meant Noether's axiomatics, then he implied, critically, that algebra resembles a mythical realm that keeps luring mathematicians with the illusory promise of great wealth; but if instead he meant what he would then come to call immanent axiomatics, then he indicated, approvingly, that algebra is a domain where genuine mathematical fortunes could actually be reached. The driving wedge between these meanings is the principle of detachment.

\section{References}

Buldt, Bernd. 2004. ``On RC 102-43-14.'' In \textit{Carnap: From Jena to L.A.},
ed. S. Awodey and C. Klein, 225--246. Chicago: Open Court.

\bigskip

Corry, Leo. 2017. ``Steht es alles wirklich schon bei Dedekind? Ideals and factorization between Dedekind and Noether.’’ In \textit{In Memoriam Richard Dedekind (1831--1916)}, ed. K. Scheel, Th. Sonar and P. Ullrich, 134--159. M\"unster: Verlag f\"ur wissenschaftliche Texte und Medien.

\bigskip

Eckes, Christophe, N. Schappacher. 2016. ``Dating the Gasthof Vollbrecht Photograph.'' Available online at https://opc.mfo.de/files/GVanglais\_def.pdf

\bigskip

Detlefsen, Michael. 2005. ``Formalism.'' In \textit{The Oxford Handbook of Philosophy of Mathematics and Logic}, ed. S. Shapiro, 236--317. Oxford: Oxford University Press.

\bigskip

Edwards, Harold M. 2005. \emph{Essays in Constructive Mathematics}. Springer.

\bigskip

Feferman, Solomon. 1964. ``Systems of Predicative Analysis.'' \textit{Journal of Symbolic Logic} 29:1--30.

\bigskip

Feferman, Solomon. 1988. ``Weyl vindicated. Das Kontinuum 70 Years later.'' In \textit{Temi e prospettive della logica e della filosofia della scienza contemporanee}, Bologna, 59--93, reprinted, with  a Postscript, in his \textit{In the Light of Logic}, 1998, 249--283. Oxford: Oxford University Press. 

\bigskip

G\"odel, Kurt. 2014. \textit{Complete Works}, vol. V. Oxford: Oxford University Press.

\bigskip

Grzegorczyk, Andrzej.  1955. ``Elementarily definable analysis.''  \textit{Fundamenta mathematicae} 41:311--338.

\bigskip

Hilbert, David. 1923. ``Die logischen Grundlagen der Mathematik.'' \textit{Mathematische Annalen} 88:151–-16. English translation in \textit{From Kant to Hilbert. A Source Book in the Foundations of Mathematics}, ed. W. B. Ewald, vol. 2, 1996, 1134–-1148. Oxford: Oxford University Press.

\bigskip

Koreuber, Mechthild. 2015. \textit{Emmy Noether, die Noether-Schule und die moderne Algebra}, Mathematik im Kontext. Springer.

\bigskip

McLarty, Colin. 2006. ``Emmy Noether's `set theoretic' topology: From Dedekind to the rise of functors.'' In \textit{The Architecture of Modern Mathematics: Essays in History and Philosophy}, ed. J. Ferreiros and J. Gray, 187--208. Oxford: Oxford University Press.

\bigskip

Moore, Gregory. 1988. ``The Emergence of First-Order Logic.'' In \textit{History and Philosophy of Modern Mathematics}, ed. W. Aspray and Ph. Kitcher, 95--135. University of Minnesota Press.

\bigskip

Noether, Emmy. 1907. ``On Complete Systems of Invariants for Ternary Biquadratic Forms.'' \textit{Journal fur die reine und angewandte Mathematik} 134:23--90.

\bigskip

Noether, Emmy. 1921. ``Ideal Theory in Rings.'' \textit{Mathematische Annalen} 83:24--66.

\bigskip

Noether, Emmy,  Schmeidler, Werner. 1920. ``On Modules in Non-Commutative Fields, Especially from Differential and Difference Expression.'' \textit{Mathematische Zeitschrift} 8:1--35.

\bigskip

Roquette, Peter. 2008. ``Emmy Noether and Hermann Weyl.'' Available at https://www.mathi.uni-heidelberg.de/~roquette/manu.html

\bigskip

Roquette, Peter. 2010. ``In Memoriam Ernst Steinitz (1871-1928).'' Available at https://www.mathi.uni-heidelberg.de/~roquette/manu.html

\bigskip

Stein, Howard. 1988. ``Logos, Logic, and Logistik\'e: Some Philosophical Remarks on Nineteenth-Century Transformation of Mathematics.'' In \textit{History and Philosophy of Modern Mathematics}, ed. W. Aspray and Ph. Kitcher, 238--259. University of Minnesota Press.

\bigskip

Steinitz, Ernst. 1910. ``Algebraische Theorie der K\"orper.'' \textit{Journal f\"ur die reine und angewandte Mathematik} 137:167--309.

\bigskip

Textor, Mark. 2019. ``Correctness First:
Brentano on Judgment and Truth.'' In \textit{The Act and Object
of Judgment
Historical and Philosophical Perspectives}, ed. B. Ball
and Ch. Schuringa, 129--150. Routledge.

\bigskip

Toader, Iulian D. 2011. \textit{Objectivity Sans Intelligibility: Hermann Weyl’s Symbolic Constructivism}.
PhD diss., University of Notre Dame.

\bigskip

Toader, Iulian D. 2013. ``Concept formation and scientific objectivity: Weyl's turn against Husserl''. In \textit{HOPOS: The
Journal of the International Society for the History of Philosophy of Science}, 3, 281--305.

\bigskip

Toader, Iulian D. 2014. ``Why Did Weyl Think that Formalism's Victory Against
Intuitionism Entails a Defeat of Pure Phenomenology?'' In \textit{History and Philosophy of Logic}, 35:2,
198--208.

\bigskip

Toader, Iulian D.
2016. ``Why Did Weyl Think that Dedekind’s Norm of Belief in Mathematics is Perverse?'' In \textit{Early Analytic Philosophy – New Perspectives on the Tradition}. The Western Ontario Series in Philosophy of Science, vol. 80, pp. 445--451.

\bigskip

van der Waerden, Bartel Leendert. 1935. ``Nachruf auf Emmy Noether.'' \textit{Mathematische Annalen} 111:469--474.

\bigskip

Weyl, Hermann. 1918. \textit{Das Kontinuum. Kritische Untersuchungen \"uber die Grundlagen der Analysis}. In \textit{Das Kontinuum und andere Monographien}, Chelsea Publishing Company, 1960. Translated as \textit{The Continuum}, 1987. Thomas Jefferson University Press.

\bigskip

Weyl, Hermann. 1927. \textit{Philosophie der Mathematik und Naturwissenschaften}. Berlin: Oldenbourg. Translated as \textit{Philosophy of Mathematics and Natural Science}. Princeton University Press. Reprinted in 2009 with an introduction by Frank Wilczek.

\bigskip

Weyl, Hermann. 1932. ``Topologie und abstrakte Algebra als zwei Wege mathematischen Verst\"{a}ndnisses.'' \textit{Unterrichtsbl\"{a}tter f\"{u}r Mathematik und Naturwissenschaften} 38:177--188. Reprinted in Weyl 1968, III, 348--358. Translated as ``Topology and Abstract Algebra as Two Roads of Mathematical Comprehension.'' \textit{American Mathematical Monthly}, 1995, 102:453-460, 102:646--651. Reprinted in Weyl 2013, 33--48.

\bigskip

Weyl, 102. 1935. ``Emmy Noether, 1882--1935.'' \textit{Scripta mathematica} 3:201--220. Reprinted in Weyl 1968, III, 425--444, and in Weyl 2013, 49--66.

\bigskip

Weyl, Hermann. 1940. ``The Mathematical Way of Thinking.'' \textit{Science} 92:437--446. Reprinted in Weyl 1968, III, 710--718, and in Weyl 2013, 67--84.

\bigskip

Weyl, Hermann. 1951. ``A Half-Century of Mathematics.'' \textit{American Mathematical Monthly} 58:523--553. Reprinted in Weyl 1968, IV, 464--494, and in Weyl 2013, 159--190.

\bigskip

Weyl, Hermann. 1953. ``Axiomatic versus Constructive Procedures in Mathematics.'' \textit{The Mathematical Intelligencer} 1985, 7/4:10--17. Reprinted in Weyl 2013, 191--202.

\bigskip

Weyl, Hermann. 1954. ``Address of the President of the Fields Medal Committee.'' In \textit{Proceedings of the International Congress of Mathematicians}. Reprinted in Weyl 1968, IV, 609--622.

\bigskip

Weyl, Hermann. 1968. \textit{Gesammelte Abhandlungen}, vol I-IV, ed. Chandrasekharan, Komaravolu. Berlin: Springer.

\bigskip

Weyl, Hermann. 2013. \textit{Levels of Infinity: Selected Writings on Mathematics and Philosophy}, ed. by Peter Pesic. Dover Publications.

\bigskip

Zilber, Boris. 2010. \textit{Zariski Geometries}.
Cambridge University Press.

\end{document}